\documentclass[11pt,reqno]{amsproc}
\usepackage[all]{xy}
\usepackage[colorlinks,urlcolor=blue,citecolor=blue,linkcolor=blue]{hyperref}
\usepackage{natbib,xcolor,graphicx,newtx,bm,anysize}
\usepackage{caption,subcaption,eqnarray}

\marginsize{15mm}{15mm}{6mm}{5mm}

\theoremstyle{theorem}
\newtheorem{theorem}{Theorem}[section]
\newtheorem{proposition}[theorem]{Proposition}

\theoremstyle{remark}

\newtheorem{remark}[theorem]{\it Remark}

\theoremstyle{definition}

\AtBeginDocument{.}

\numberwithin{equation}{section}
\allowdisplaybreaks

\newcommand{\assign}{:=}
\newcommand{\cdummy}{\cdot}

\newcommand{\mathe}{\mathrm{e}}

\newcommand{\tmop}[1]{\ensuremath{\operatorname{#1}}}

\newenvironment{enumeratenumeric}{\begin{enumerate} }{\end{enumerate}}

\begin{document}
\title[Multiple solutions for Schr{\"o}dinger--Poisson--Slater equations]{Multiple solutions for Schr{\"o}dinger--Poisson--Slater equations with
critical growth}
\author{Shibo Liu\vspace{-1em}}
\dedicatory{Department of Mathematical Sciences, Florida Institute of Technology\\
Melbourne, FL 32901, USA\\[1em]
\sf Dedicated to Professor Shujie Li on the Occasion of his 85th Birthday}
\thanks{Emails: \texttt{\bfseries sliu@fit.edu} (S. Liu)}
\begin{abstract}
We obtain multiple solutions for the zero mass
  Schr{\"o}dinger--Poisson--Slater equation
  \[ - \Delta u + \left( \frac{1}{4 \pi | x |} \ast u^2 \right) u = \lambda g
     (x) | u |^{p - 2} u + | u |^{6 - 2} u \text{, \qquad} u \in
     \mathcal{D}^{1, 2} (\mathbb{R}^3) \text{} \]
  for $\lambda \gg 1$, where $p \in (4, 6)$ and $g \in L^{6 / (6 - p)}
  (\mathbb{R}^3)$. The crucial $(\tmop{PS})_c$ condition is verified using a
  simpler method. Similar multiplicity result is also obtained for related
  equation with an external potential.
\end{abstract}
\maketitle

\section{Introduction}

We consider the following zero mass Schr{\"o}dinger--Poisson--Slater equation
\begin{equation}
  - \Delta u + \left( \frac{1}{4 \pi | x |} \ast u^2 \right) u = \lambda g (x)
  | u |^{p - 2} u + | u |^{6 - 2} u \text{, \qquad} u \in \mathcal{D}^{1, 2}
  (\mathbb{R}^3) \text{.} \label{e1}
\end{equation}
Note that the exponent $6$ in the last term is the critical exponent for the
Sobolev embedding.

Nonlocal elliptic equations like (\ref{e1}) and its counterpart (\ref{e11})
arise from finding standing waves $\psi (t, x) = \mathe^{- i \omega t} u
(x)$ for the following nonlocal Schr{\"o}dinger equation
\[ i \partial_t \psi = - \Delta \psi + U (x) \psi + \left( \frac{1}{4 \pi | x
   |} \ast | \psi |^2 \right) \psi - f (x, | \psi |) \psi \text{, \quad} (t,
   x) \in \mathbb{R}_+ \times \mathbb{R}^3 \text{,} \]
which comes from \ an approximation of the Hartree--Fock model of a quantum
many-body system of electrons, in which $| \psi |^2$ is the density of
electrons and the nonlocal convolution term represents the Coulombic repulsion
between the electrons.

{\cite{MR2679375}} and {\cite{MR2902293}} studied
(\ref{e1}) for the case that the nonlinearity (the right hand side) is a pure
power function $| u |^{q - 2} u$ with $q \in (2, 6)$. 
\cite{MR3912770} studied (\ref{e1}) for the case that $g (x) \equiv 1$, they
obtained ground state solution for $p \in (3, 6)$ and positive radial solution
for $p \in \left( \frac{18}{7}, 3 \right)$, see also {\cite{Gu2024}}.

More recently, for $g (x) \equiv 1$ and $p \in (3, 6)$, 
{\citet[Theorem 1.28]{mercuri2024variationalmethodsscaledfunctionals}} obtained
arbitrarily many solutions for the equation (\ref{e1}), provided $\lambda$ is
large enough. Since $g (x)$ is constant, they can work with the radial
subspace $E_r (\mathbb{R}^3)$ of the Coulomb--Sobolev space $E
(\mathbb{R}^3)$. The variational functional enjoys nice $3$-scaling property,
which enable them to apply their critical point theorem for scaled functionals
{\cite[Corollary 2.34]{mercuri2024variationalmethodsscaledfunctionals}}.

In this paper we study the case that $g$ is not a constant, not even radially
symmetric. Thus, we have to work on the general Coulomb--Sobolev space $E
(\mathbb{R}^3)$ and the crucial $3$-scaling property is lost. Nevertheless, we
still obtain similar multiplicity result with narrower range of $p$. We assume
that $g$ satisfies
\begin{enumerate}
  \item[$(g)$] there is $p \in (4, 6)$ such that $g \in L^{6 / (6 - p)}
  (\mathbb{R}^3)$, $g \geqslant 0$, $\Omega \assign \{ g > 0 \}$ is nonempty
  open subset of $\mathbb{R}^3$.
\end{enumerate}
Then, we have the following theorem.

\begin{theorem}
  \label{t1}Let $g$ satisfy the condition $(g)$. Given $m \in \mathbb{N}$,
  there is $\lambda_m > 0$ such that \eqref{e1} has $m$ pairs of solutions
  with positive energy for all $\lambda \geqslant \lambda_m$.
\end{theorem}

Our proof of this theorem is based on a critical point theorem of 
{\citet[Theorem 2.1]{perera2024abstractmultiplicitytheoremsapplications}}. Like
almost all critical point theorem, Palais--Smale $(\tmop{PS})$ condition is
crucial for applying this theorem. Since (\ref{e1}) is of critical growth, the
most we can expect is local $(\tmop{PS})$ condition, that is $(\tmop{PS})_c$
for all $c \in (0, c^{\ast})$ for some $c^{\ast} > 0$. In  
{\cite{mercuri2024variationalmethodsscaledfunctionals}} the proof of
$(\tmop{PS})_c$ for the case $g\equiv1$ depends on the Pohozaev identity {\cite[Lemma
2.4]{MR4292779}}. For our case that $g$ is not constant, we will give a
simpler proof in Section \ref{s1}. In Section \ref{s} we first recall the
Coulomb--Sobolev space $E (\mathbb{R}^3)$ introduced by 
{\cite{MR2679375}}, then present the proof of Theorem \ref{t1}. In Section
\ref{s3}, we present similar result for Schr{\"o}dinger--Poisson--Slater
equation with an external potential (see Eq. (\ref{e11})). Finally, in Section
\ref{sss} we present some variants of the results we have obtained so far.

\section{Variational setting and proof of Theorem \ref{t1}}\label{s}

Instead of the standard Sobolev space $H^1 (\mathbb{R}^3)$, the correct
functional space for studying the zero mass problem (\ref{e1}) is the
Coulomb--Sobolev space $E (\mathbb{R}^3)$ introduced by 
{\cite{MR2679375}}, where $E (\mathbb{R}^3)$ is the vector space
\[ E = E (\mathbb{R}^3) = \left\{ u \in \mathcal{D}^{1, 2} (\mathbb{R}^3)
   \left| \iint \frac{u^2 (x) u^2 (y)}{| x - y |} < \infty \right. \right\}
   \text{} \]
equipped with the norm
\[ \| u \| = \left[ \int | \nabla u |^2 + \left( \iint \frac{u^2 (x) u^2
   (y)}{| x - y |} \right)^{1 / 2} \right]^{1 / 2} \text{.} \]
Here and in what follows, unless stated explicitly, all integrals are taken
over $\mathbb{R}^3$, all double integrals are taken with respect to $(x, y)$ over
$\mathbb{R}^3 \times \mathbb{R}^3$.

It has been proved in {\cite[Theorem 1.5]{MR2679375}} that $(E, \|
\cdot \|)$ is a uniformly convex Banach space which is embedded in $L^q (\mathbb{R}^3)$
continuously for $q \in \left[3, 6 \right]$.

We consider the functional $\Phi : E \rightarrow \mathbb{R}$,
\[ \Phi (u) = \frac{1}{2} \int | \nabla u |^2 + \frac{1}{16 \pi} \iint
   \frac{u^2 (x) u^2 (y)}{| x - y |} - \frac{\lambda}{p} \int g | u |^p -
   \frac{1}{6} \int | u |^6 \text{.} \]
Then, it is well known that $\Phi \in C^1 (E)$ with derivative given by
\[ \langle \Phi' (u), v \rangle = \int \nabla u \cdummy \nabla v + \frac{1}{4
   \pi} \iint \frac{u^2 (x) u (y) v (y)}{| x - y |} - \frac{\lambda}{p} \int g
   | u |^{p - 2} u v - \int | u |^{6 - 2} u v \text{.} \]
Hence, critical points of $\Phi$ are weak solutions of the problem (\ref{e1}).
By regularity results, weak solutions are classical solutions.

Therefore, we will focus on finding multiple critical points of $\Phi$. For
this purpose, we need the following critical point theorem of
{\cite{perera2024abstractmultiplicitytheoremsapplications}}. For a symmetric
subset $A$ of $E \backslash \{ 0 \}$, we denote by $i (A)$ the cohomological
index of $A$, which was introduced by {\cite{MR0478189}}.
If $A$ is homeomorphic to the unit sphere $S^{m - 1}$ in $\mathbb{R}^m$, then
$i (A) = m$.

\begin{proposition}[{\citet[Theorem
2.1]{perera2024abstractmultiplicitytheoremsapplications}}]
  \label{p1} Let $E$ be a Banach space, $\Phi : E \rightarrow \mathbb{R}$ be
  an even $C^1$-functional satisfying $(\tmop{PS})_c$ for $c \in (0,
  c^{\ast})$ being $c^{\ast}$ some positive constant. If $0$ is a strict local
  minimizer of $\Phi$ and there are $R > 0$ and a compact symmetric set $A
  \subset \partial \mathfrak{B}_R$, where $\mathfrak{B}_R$ is the $R$-ball in
  $E$, such that $i (A) = m$,
  \begin{equation}
    \max_A \Phi \leqslant 0 \text{, \qquad} \max_B \Phi < c^{\ast} \text{,}
    \label{3e}
  \end{equation}
  where $B = \{ t u \mid t \in [0, 1], u \in A \}$, then $\Phi$ has $m$ pairs
  of nonzero critical points with positive critical values.
\end{proposition}

\subsubsection*{Proof of Theorem \ref{t1}}As has been pointed out by
{\cite{MR2902293}}, it is clear that
\[ \frac{1}{2} \| u \|^4 \leqslant \int | \nabla u |^2 + \iint \frac{u^2 (x)
   u^2 (y)}{| x - y |} \text{\qquad if } \| u \| \leqslant 1 \text{.} \]
Since $p > 4$, using the above inequality it is clear that $u = 0$ is a strict
local minimizer of $\Phi $.

Let
\[ S = \inf_{u \in \mathcal{D}^{1, 2} (\mathbb{R}^3) \backslash \{ 0 \}}
   \frac{| \nabla u |_2^2}{| u |_6^2} \text{, \qquad} c^{\ast} = \frac{1}{3}
   S^{3 / 2} \text{.} \]
    In Section \ref{s1} we will show that $\Phi$
satisfies $(\tmop{PS})_c$ for $c \in (0, c^{\ast})$. To conclude the proof of
Theorem \ref{t1}, it suffices to find the subsets $A$ and $B$ satisfying the
geometric assumption (\ref{3e}) for any given $m \in \mathbb{N}.$ We will
adapt the argument used in {\cite{MR4785363}}, where a $(p, q)$-Laplacian
equation
\[ \left\{ \begin{array}{l}
     - \Delta_p u - \Delta_q u = \lambda h (x) | u |^{r - 2} u + g (x) | u
     |^{p^{\ast} - 2} u \text{,}\\
     u \in \mathcal{D}^{1, p} (\mathbb{R}^N) \cap \mathcal{D}^{1, q}
     (\mathbb{R}^N)
   \end{array} \right. \]
is considered.

Given $m \in \mathbb{N}$, let
\[ Z = \{ u \in E \mid \tmop{supp} u \subset \Omega \} \text{,} \]
$Z_m$ be an $m$-dimensional subspace of $Z$. Since $g > 0$ on $\Omega$,
\[ [u]_g = \left( \frac{1}{p} \int g | u |^p \right)^{1 / p} \text{} \]
is a norm on $Z$ and $Z_m$. For $u \in Z_m$ we have
\begin{align}
  \Phi (u) & \leqslant  \frac{1}{2} \| u \|^2 + \frac{1}{16 \pi} \| u \|^4 -
  \lambda [u]_g^p - \frac{1}{6} | u |_6^6 \nonumber\\
  & \leqslant  \frac{1}{2} \| u \|^2 + \frac{1}{16 \pi} \| u \|^4 - \lambda
  a_1 \| u \|^p - a_2 \| u \|^6 \text{} \label{e2} 
\end{align}
because all norms on $Z_m$ are equivalent. Take $R > 1$ such that
\begin{equation}
  f (R) : = \frac{1}{2} R^2 + \frac{1}{16 \pi} R^4 - a_2 R^6 < 0 \text{.}
  \label{e3}
\end{equation}
Let $A = Z_m \cap \partial \mathfrak{B}_R$, then $i (A) = m$. If $\lambda >
0$, then for any $u \in A$, from (\ref{e2}) we have $\Phi (u) \leqslant f
(R)$. Thus
\[ \max_A \Phi < 0 \text{.} \]
For the function $f$ defined in (\ref{e3}), there is $\delta \in (0, R)$ such
that $f (s) < c^{\ast}$ for all $s \in [0, \delta]$. Set
\[ \lambda_m = 1 + \max_{s \in [\delta, R]} \left| \frac{f (s) - c^{\ast}}{a_1
   s^p} \right| \text{.} \]
Then if $\lambda \geqslant \lambda_m$ we have
\[ f (s) - \lambda a_1 s^p < c^{\ast} \text{\qquad for } s \in [\delta, R]
   \text{.} \]
Therefore, for $u \in A$,
\begin{enumeratenumeric}
  \item if $t \in \left[ \frac{\delta}{R}, 1 \right]$, then $\| t u \| \in
  [\delta, R]$,
  \[ \Phi (t u) \leqslant f (\| t u \|) - \lambda a_1 \| t u \|^p < c^{\ast}
     \text{;} \]
  \item if $t \in \left[ 0, \frac{\delta}{R} \right]$, then $\| t u \|
  \leqslant \delta$ and $\Phi (t u) \leqslant f (\| t u \|) < c^{\ast}$.
\end{enumeratenumeric}
From this, we deduce that for $B = \{ t u \mid t \in [0, 1], u \in A \}$ there
holds
\[ \max_B \Phi < c^{\ast} \text{.} \]
By Proposition \ref{p1}, $\Phi$ has $m$ pairs of nonzero critical points, and
(\ref{e1}) has $m$-pairs of nontrivial solutions.

\section{$(\tmop{PS})_c$ condition}\label{s1}

In this section we show that $\Phi$ satisfies $(\tmop{PS})_c$ condition for
all $c \in (0, c^{\ast})$. Let $\{u_n\}$ be a $(\tmop{PS})_c$ sequence with $c \in
(0, c^{\ast})$, that is
\begin{align}
  \Phi (u_n) & =  \frac{1}{2} \int | \nabla u_n |^2 + \frac{1}{16 \pi} \iint
  \frac{u_n^2 (x) u_n^2 (y)}{| x - y |} - \frac{\lambda}{p} \int g | u_n |^p -
  \frac{1}{6} \int | u_n |^6 \rightarrow c \text{,} \nonumber\\
  \langle \Phi' (u_n), v \rangle & =  \int \nabla u_n \cdummy \nabla v +
  \frac{1}{4 \pi} \iint \frac{u_n^2 (x) u_n (y) v (y)}{| x - y |} \nonumber\\
  &   \hspace{4em} - \lambda \int g | u_n |^{p - 2} u_n v - \int | u_n |^{6
  - 2} u_n v = o (\| v \|) \text{.} \label{e8} 
\end{align}
Since $p \in (4, 6)$, we may take $\mu \in (4, p)$. Then for $n \gg 1$ we have 
\begin{align}
  c + 1 & \geqslant  \Phi (u_n) - \frac{1}{\mu} \langle \Phi' (u_n), u_n
  \rangle \nonumber\\
  & =  \left( \frac{1}{2} - \frac{1}{\mu} \right) \int | \nabla u_n |^2 +
  \left( \frac{1}{16 \pi} - \frac{1}{4 \mu \pi} \right) \iint \frac{u_n^2 (x)
  u_n^2 (y)}{| x - y |} \nonumber\\
  &   \hspace{4em} + \left( \frac{\lambda}{\mu} - \frac{\lambda}{p} \right)
  \int g | u_n |^p + \left( \frac{1}{\mu} - \frac{1}{6} \right) \int | u_n |^6
  \nonumber\\
  & \geqslant  \left( \frac{1}{2} - \frac{1}{\mu} \right) \int | \nabla u_n
  |^2 + \left( \frac{1}{16 \pi} - \frac{1}{4 \mu \pi} \right) \iint
  \frac{u_n^2 (x) u_n^2 (y)}{| x - y |} \text{.} \label{bdd} 
\end{align}
Since $\mu > 4$, the coefficients of the integrals at the end are positive. It
follows that $\{u_n\}$ is bounded.

Up to a subsequence we may assume $u_n \rightharpoonup u$ in $E$, and
\[ u_n \rightharpoonup u \text{\quad in } \mathcal{D}^{1, 2} (\mathbb{R}^3)
   \text{, \qquad} u_n \rightarrow u \text{\quad a.e. on } \mathbb{R}^3
   \text{.} \]
Since $\left\{| u_n |^{6 - 2} u_n\right\}$ is bounded in $L^{6 / 5} (\mathbb{R}^3)$ and $|
u_n |^{6 - 2} u_n \rightarrow | u |^{6 - 2} u$ a.e.\ on $\mathbb{R}^3$, using
{\cite[Page 487]{MR699419}} we have $| u_n |^{6 - 2} u_n \rightharpoonup | u
|^{6 - 2} u$ in $L^{6 / 5} (\mathbb{R}^3)$. Therefore
\begin{equation}
  \int | u_n |^{6 - 2} u_n v \rightarrow \int | u |^{6 - 2} u v \text{\qquad
  for } v \in E \text{.} \label{A1}
\end{equation}
Similarly, since $\left\{| u_n |^p\right\}$ and $\left\{| u_n |^{p - 2} u_n v\right\}$ are bounded in $L^{6
/ p} (\mathbb{R}^3)$, using {\cite[Page 487]{MR699419}} again we have
\[
| u_n |^p
\rightharpoonup | u |^p\text{\quad in }L^{6 / p} (\mathbb{R}^3)\text{,\qquad}| u_n |^{p - 2} u_n v
\rightharpoonup | u |^{p - 2} u v\text{\quad in }L^{6 / p} (\mathbb{R}^3)\text{.}
\]
Because $g \in
L^{6 / (6 - p)} (\mathbb{R}^3)$, the dual of $L^{6/p}(\mathbb{R}^3)$ we deduce
\begin{equation}
  \int g | u_n |^p \rightarrow \int g | u |^p \text{, \qquad} \int g | u_n
  |^{p - 2} u_n v \rightarrow \int g | u |^{p - 2} u v \text{.} \label{A2}
\end{equation}
Moreover, since $u_n \rightharpoonup u$ in $E$, using {\cite[Lemma 2.3]{MR2902293}} we have
\begin{equation}
  \iint \frac{u_n^2 (x) u_n (y) v (y)}{| x - y |} \rightarrow \iint \frac{u^2
  (x) u (y) v (y)}{| x - y |} \text{.} \label{A3}
\end{equation}
Using (\ref{A1}), (\ref{A2}) and (\ref{A3}) we deduce
\begin{align}
  0 & =  \lim_{n \rightarrow \infty} \langle \Phi' (u_n), v \rangle
  \nonumber\\
  & =  \int \nabla u \cdummy \nabla v + \frac{1}{4 \pi} \iint \frac{u^2 (x)
  u (y) v (y)}{| x - y |} - \lambda \int g | u |^{p - 2} u v - \int | u |^{6 -
  2} u v \label{e4} \\
  & =  \langle \Phi' (u), v \rangle \text{.} \nonumber
\end{align}
So $\Phi' (u) = 0$. We also have $\Phi (u) \geqslant 0$ because
\begin{align*}
  4 \Phi (u) & =  4 \Phi (u) - \langle \Phi' (u), u \rangle\\
  & =  \int | \nabla u |^2 + \lambda \left( 1 - \frac{4}{p} \right) \int g |
  u |^p + \left( 1 - \frac{4}{6} \right) \int | u |^6 \text{.}
\end{align*}

Let $v_n = u_n - u$. By Brezis--Lieb lemma {\cite[Theorem 1]{MR699419}} (see
also {\cite[Lemma 1.32]{MR1400007}},
\begin{equation}
  \int | u_n |^6 = \int | u |^6 + \int | v_n |^6 + o (1) \text{.} \label{e5}
\end{equation}
Using this and
\begin{align}
  \int | \nabla u_n |^2 & =  \int | \nabla u |^2 + \int | \nabla v_n |^2 + o
  (1) \text{,} \label{e6} \\
  \int g | u_n |^p & =  \int g | u |^p + o (1) \text{,} \label{e7} 
\end{align}
as well as $\Phi (u_n) \rightarrow c$, we deduce
\begin{align}
   4 \Phi (u)& + 2 \int | \nabla v_n |^2 - \frac{2}{3} \int | v_n |^6\nonumber\\
  & +
  \frac{1}{4 \pi} \iint \frac{u_n^2 (x) u_n^2 (y)}{| x - y |} - \frac{1}{4
  \pi} \iint \frac{u^2 (x) u^2 (y)}{| x - y |} \rightarrow 4 c \text{.}
  \label{Q0} 
\end{align}
On the other hand, it follows from $\langle \Phi' (u_n), u_n \rangle
\rightarrow 0$ and $\Phi' (u) = 0$ that
\begin{align}
  0 & =  \langle \Phi' (u_n), u_n \rangle + o (1) \nonumber\\
  & =  \int | \nabla u_n |^2 + \frac{1}{4 \pi} \iint \frac{u_n^2 (x) u_n^2
  (y)}{| x - y |} - \lambda \int g | u_n |^p - \int | u_n |^6 + o (1)
  \nonumber\\
  & =  \int | \nabla u |^2 + \int | \nabla v_n |^2 + \frac{1}{4 \pi} \iint
  \frac{u_n^2 (x) u_n^2 (y)}{| x - y |} \nonumber\\
  &   \hspace{4em} - \lambda \int g | u |^p - \int | u |^6 - \int | v_n |^6
  + o (1) \nonumber\\
  & =  \langle \Phi' (u), u \rangle + \int | \nabla v_n |^2 - \int | v_n |^6
  \nonumber\\
  &   \hspace{4em} + \frac{1}{4 \pi} \iint \frac{u_n^2 (x) u_n^2 (y)}{| x -
  y |} - \frac{1}{4 \pi} \iint \frac{u^2 (x) u^2 (y)}{| x - y |} + o (1)
  \nonumber\\
  & =  \int | \nabla v_n |^2 - \int | v_n |^6 + \frac{1}{4 \pi} \iint
  \frac{u_n^2 (x) u_n^2 (y)}{| x - y |} - \frac{1}{4 \pi} \iint \frac{u^2 (x)
  u^2 (y)}{| x - y |} + o (1) \text{.} \label{P1} 
\end{align}
Using (\ref{Q0}), (\ref{P1}) we get (the double integral terms cancel)
\begin{equation}
  4 \Phi (u) + | \nabla v_n |_2^2 + \frac{1}{3} | v_n |_6^6 \rightarrow 4 c
  \text{.} \label{W}
\end{equation}
Since $\left\{| \nabla v_n |_2^2\right\}$ and $\left\{| v_n |_6^6\right\}$ are bounded, up to a subsequence
we may assume that
\[ | \nabla v_n |_2^2 \rightarrow a \text{, \qquad} | v_n |_6^6 \rightarrow b
   \text{.} \]
Since $u_n \rightarrow u$ a.e.\ on $\mathbb{R}^3$, by Fatou lemma we have
\begin{equation}
  \varliminf_{n \rightarrow \infty} \iint \frac{u_n^2 (x) u_n^2 (y)}{| x - y |}
  \geqslant \iint \frac{u^2 (x) u^2 (y)}{| x - y |} \text{.}
\end{equation}
From this and (\ref{P1}) we see that $b \geqslant a$.

By the definition of $S$, $| \nabla v_n |_2^2 \geqslant S | v_n |_6^2$, which
implies
\[ a \geqslant S b^{2 / 6} \geqslant S a^{2 / 6} \text{.} \]
So either $a = 0$ or $a \geqslant S^{3 / 2}$. If $a \geqslant S^{3 / 2}$, we
deduce from (\ref{W}), $\Phi (u) \geqslant 0$ and $b \geqslant a$ that
\[ c \geqslant \frac{1}{4} a + \frac{1}{12} b \geqslant \frac{1}{3} a
   \geqslant \frac{1}{3} S^{3 / 2} = c^{\ast} \text{,} \]
which contradicts our assumption that $c < c^{\ast}$. Therefore, $a = 0$. That
is, $v_n \rightarrow 0$ in $\mathcal{D}^{1, 2} (\mathbb{R}^3)$.

Now we follow the argument in {\cite[Page
58]{mercuri2024variationalmethodsscaledfunctionals}}. By the nonlocal
Brezis--Lieb lemma {\cite[Proposition 4.1]{Mercuri2016}}, we have
\begin{equation}
  \iint \frac{u_n^2 (x) u_n^2 (y)}{| x - y |} \geqslant \iint \frac{u^2 (x)
  u^2 (y)}{| x - y |} + \iint \frac{v_n^2 (x) v_n^2 (y)}{| x - y |} + o (1)
  \text{.} \label{e9}
\end{equation}
Subtracting (\ref{e4}) with $v = u$, that is
\[ \int | \nabla u |^2 + \frac{1}{4 \pi} \iint \frac{u^2 (x) u^2 (y)}{| x - y
   |} - \lambda \int g | u |^p - \int | u |^6 = 0 \text{,} \]
from (\ref{e8}) with $v = u_n$, that is
\[ \int | \nabla u_n |^2 + \frac{1}{4 \pi} \iint \frac{u_n^2 (x) u^2_n (y)}{|
   x - y |} - \lambda \int g | u_n |^p - \int | u_n |^6 = \langle \Phi' (u_n),
   u_n \rangle = o (1) \text{,} \]
then using (\ref{e5}), (\ref{e6}), (\ref{e7}) and (\ref{e9}), as well as $v_n
\rightarrow 0$ in $\mathcal{D}^{1, 2} (\mathbb{R}^3)$, we get
\begin{align}
  \int | \nabla v_n |^2 + \iint \frac{v_n^2 (x) v_n^2 (y)}{| x - y |} &
  \leqslant  \int | v_n |^6 + o (1) \nonumber\\
  & \leqslant  S^{- 3} \left( \int | \nabla v_n |^2 \right)^3 + o (1)
  \rightarrow 0 \text{.} \label{19} 
\end{align}
This means $v_n \rightarrow 0$ in $E$, that is $u_n \rightarrow u$ in $E$.

\begin{remark}
  Unlike the proof of {\cite[Lemma 3.6]{mercuri2024variationalmethodsscaledfunctionals}} for the case $g
  \equiv 1$, our proof here does not depend on the Pohozaev identity
  {\cite[Lemma 2.4]{MR4292779}}, therefore is somewhat simpler. If $g\equiv1$ then as in \cite{mercuri2024variationalmethodsscaledfunctionals} instead of $E$ we should work on $E_r$. Thanks to the compact embedding $E_r\hookrightarrow  L^p(\mathbb{R}^3)$, (\ref{A2}) is still valid and our argument works as well.
\end{remark}

\section{Potential case}\label{s3}

The argument above can be applied to similar equations with an external
potential
\begin{equation}
  - \Delta u + V (x) u + \left( \frac{1}{4 \pi | x |} \ast u^2 \right) u =
  \lambda g (x) | u |^{p - 2} u + | u |^{q - 2} u \text{, \qquad} u \in H^1
  (\mathbb{R}^3) \text{,} \label{e11}
\end{equation}
where $4<p<q\le6$, $V : \mathbb{R}^3 \rightarrow (0, \infty)$ is any external potential
such that
\[ \| u \|_V = \left( \int (| \nabla u |^2 + V u^2) \right)^{1 / 2} \]
is equivalent to the standard $H^1$-norm. The equation (\ref{e11}) is
equivalent to the nonlinear Schr{\"o}dinger-Poisson systems with critical or subcritical
growth:
\[ \left\{ \begin{array}{ll}
     - \Delta u + V (x) u + \phi u = \lambda g (x) | u |^{p - 2} u + | u |^{q
     - 2} u & \text{in } \mathbb{R}^3 \text{,}\\
     - \Delta \phi = u^2 & \text{in } \mathbb{R}^3 \text{,}
   \end{array} \right. \]
for which, there are many results, see \cite{Zhao2009,Yao2019,Furtado2023,Kang2023,
MR4653880,Zhang2025a}. We have the following multiplicity result.

\begin{theorem}
  \label{t2}Let $g$ satisfy the condition $(g)$. Given $m \in \mathbb{N}$,
  there is $\lambda_m > 0$ such that \eqref{e11} has $m$ pairs of solutions
  with positive energy for all $\lambda \geqslant \lambda_m$.
\end{theorem}

For the proof, let $X$ be the completion of $C_0^{\infty} (\mathbb{R}^3)$ with
respect to the norm
\begin{equation}
  \| u \| = \left[ \int (| \nabla u |^2 + V u^2) + \left( \iint \frac{u^2 (x)
  u^2 (y)}{| x - y |} \right)^{1 / 2} \right]^{1 / 2} \text{.} \label{21}
\end{equation}
Clearly we have continuous embeddings $X \hookrightarrow E$ and $X
\hookrightarrow H^1 (\mathbb{R}^3)$. Therefore we can define $\Phi : X \rightarrow \mathbb{R}$ via
\[ \Phi (u) = \frac{1}{2} \int (| \nabla u |^2 + V u^2) + \frac{1}{16 \pi}
   \iint \frac{u^2 (x) u^2 (y)}{| x - y |} - \frac{\lambda}{p} \int g | u |^p
   - \frac{1}{q} \int | u |^q \text{.} \]
Then $\Phi \in C^1 (X)$ and critical points of $\Phi$ are solutions of
(\ref{e11}). Let
\begin{equation}
  S_q = \inf_{u \in H^1 (\mathbb{R}^3) \backslash \{ 0 \}} \frac{\| u
  \|_V^2}{| u |_q^2} \text{, \qquad} c^{\ast} = \frac{1}{3} S_q^{3 / 2}
  \text{.} \label{S11}
\end{equation}
Similar to the proof of Theorem \ref{t1}, for $m \in \mathbb{N}$ let
\[ Z = \{ u \in X \mid \tmop{supp} u \subset \Omega \} \]
and $Z_m$ be an $m$-dimensional subspace of $Z$. As in the proof of Theorem
\ref{t1}, there is $R > 0$ such that (\ref{3e}) holds for $A = Z_m \cap
\partial \mathfrak{B}_R$ and
\[ B = \{ t u \mid t \in [0, 1], u \in A \} \text{.} \]
Since $i (A) = m$, we see that the geometric conditions of Proposition
\ref{p1} hold. To get $m$-pairs of critical points for $\Phi$, it suffices to
verify $(\tmop{PS})_c$ for $c \in (0, c^{\ast})$ with $c^*$ now given in (\ref{S11}).

Thus, let $\{u_n\}$ be a $(\tmop{PS})_c$ sequence with $c \in (0, c^{\ast})$.
Similar to (\ref{bdd}), we have
\begin{align*}
  c + 1  \geqslant  \left( \frac{1}{2} - \frac{1}{\mu} \right) \| u \|_V^2 +
  \left( \frac{1}{16 \pi} - \frac{1}{4 \mu \pi} \right) \iint \frac{u_n^2 (x)
  u_n^2 (y)}{| x - y |} \text{.}
\end{align*}
It follows that $u_n$ is bounded in $X$.

Thanks to the continuous embeddings $X \hookrightarrow E$ and $X
\hookrightarrow H^1 (\mathbb{R}^3)$, up to a subsequence we have
\[ u_n \rightharpoonup u \text{\quad in } E \text{, \qquad and\qquad} u_n
   \rightharpoonup u \text{\quad in } H^1 (\mathbb{R}^3) \text{.} \]
Hence, replacing the exponent $6$ by $q$, the argument between (\ref{A1}) through (\ref{W})  remains valid except
the $\mathcal{D}^{1, 2}$-norm needs to be replaced by $\| \cdummy \|_V$, the
equivalent $H^1$-norm. In particular, (\ref{W}) now reads
\begin{equation}
  4 \Phi (u) + \| v_n \|_V^2 + \frac{1}{3} | v_n |_q^q \rightarrow 4 c
  \text{,} \label{q7}
\end{equation}
being $v_n = u_n - u$. Note that $u$ is a critical point of $\Phi : X
\rightarrow \mathbb{R}$ with $\Phi (u) \geqslant 0$.

As before, assuming
\[ \| v_n \|_V^2 \rightarrow a \text{, \qquad} | v_n |_q^q \rightarrow b
   \text{.} \]
Then $b \geqslant a$. Combining $\| v_n \|_V^2 \geqslant S_q | v_n |_q^2$, a
consequence of ($\ref{S11}$), we get $a \geqslant S_q a^{1 / 3}$. If $a
\geqslant S_q^{3 / 2}$ we get from (\ref{q7}) and $b \geqslant a$ that
\[ c \geqslant \frac{1}{4} a + \frac{1}{12} b \geqslant \frac{1}{3} a
   \geqslant \frac{1}{3} S_q^{3 / 2} = c^{\ast} \text{,} \]
contradicting $c \in (0, c^{\ast})$. Thus $a = 0$ and $v_n \rightarrow 0$ in
$H^1 (\mathbb{R}^3)$. The estimate (\ref{19}) now reads
\[ \| v_n \|_V^2 + \iint \frac{v_n^2 (x) v_n^2 (y)}{| x - y |} \leqslant
   S_q^{- 3} \| v_n \|_V^3 + o (1) \rightarrow 0 \text{,} \]
which means $u_n \rightarrow u$ in $X$.

\section{Variants of Theorem \ref{t2}}\label{sss}

Checking the proofs of Theorems \ref{t1} and \ref{t2}, we see that the
condition that $g \in L^{6 / (6 - p)} (\mathbb{R}^3)$ is only used to ensure
\begin{equation}
  \int g | u_n |^p \rightarrow \int g | u |^p \text{, \qquad} \int g | u_n
  |^{p - 2} u_n v \rightarrow \int g | u |^{p - 2} u v \text{} \label{23}
\end{equation}
for $u_n \rightharpoonup u$ in $E$ or $X$, see (\ref{A2}). Therefore, we can replace
this conditions by other conditions ensuring (\ref{23}). For example, it is
well known that (\ref{23}) holds provided
\begin{equation}
  \lim_{| x | \rightarrow \infty} g (x) = 0 \text{.} \label{24}
\end{equation}
So we have the following variant of Theorems \ref{t1} and \ref{t2}.

\begin{theorem}
  Assume that the continuous function $g : \mathbb{R}^3 \rightarrow (0,
  \infty)$ satisfy \eqref{24}, $p \in (4, 6)$. Given $m \in \mathbb{N}$, there
  is $\lambda_m > 0$ such that both \eqref{e1} and \eqref{e11} have $m$ pairs
  of solutions with positive energy for all $\lambda \geqslant \lambda_m$.
\end{theorem}

On the other hand, if the potential $V$ is coercive
\begin{equation}
  \lim_{| x | \rightarrow \infty} V (x) = + \infty \text{,} \label{25}
\end{equation}
then by {\cite{MR1349229}} we have a compact embedding $H_V
\hookrightarrow L^2 (\mathbb{R}^3)$. As a consequence the embedding $X
\hookrightarrow L^2 (\mathbb{R}^3)$ is also compact and (\ref{23}) is valid
provided $g \in L^{\infty} (\mathbb{R}^3)$. Hence, for $V$ satisfying
(\ref{25}) the same multiplicity result is true assuming $g \in L^{\infty}
(\mathbb{R}^3)$ and $p \in (4, 6)$.


\begin{thebibliography}{19}
\providecommand{\natexlab}[1]{#1}

\bibitem[{Bartsch \& Wang(1995)}]{MR1349229}
Bartsch T, Wang ZQ (1995).
\newblock \href{http://dx.doi.org/10.1080/03605309508821149}{\emph{Existence
  and multiplicity results for some superlinear elliptic problems on {${\bf
  R}^N$}}}.
\newblock Comm. Partial Differential Equations, 20(9-10) 1725--1741.

\bibitem[{Br{\'e}zis \& Lieb(1983)}]{MR699419}
Br{\'e}zis H, Lieb E (1983).
\newblock \href{http://dx.doi.org/10.2307/2044999}{\emph{A relation between
  pointwise convergence of functions and convergence of functionals}}.
\newblock Proc. Amer. Math. Soc., 88(3) 486--490.

\bibitem[{Dutko \emph{et~al.}(2021)Dutko, Mercuri \& Tyler}]{MR4292779}
Dutko T, Mercuri C, Tyler TM (2021).
\newblock \href{https://doi.org/10.1007/s00526-021-02045-y}{\emph{Groundstates
  and infinitely many high energy solutions to a class of nonlinear
  {S}chr\"odinger-{P}oisson systems}}.
\newblock Calc. Var. Partial Differential Equations, 60(5) Paper No. 174, 46.

\bibitem[{Fadell \& Rabinowitz(1978)}]{MR0478189}
Fadell ER, Rabinowitz PH (1978).
\newblock \emph{Generalized cohomological index theories for {L}ie group
  actions with an application to bifurcation questions for {H}amiltonian
  systems}.
\newblock Invent. Math., 45(2) 139--174.

\bibitem[{Furtado \emph{et~al.}(2023)Furtado, Wang \& Zhang}]{Furtado2023}
Furtado MF, Wang Y, Zhang Z (2023).
\newblock \href{https://doi.org/10.1016/j.jmaa.2023.127252}{\emph{Positive and
  nodal ground state solutions for a critical {S}chr\"{o}dinger-{P}oisson
  system with indefinite potentials}}.
\newblock J. Math. Anal. Appl., 526(2) Paper No. 127252, 23.

\bibitem[{Gu \& Liao(2024)}]{Gu2024}
Gu Y, Liao F (2024).
\newblock \href{https://doi.org/10.1007/s12220-024-01656-z}{\emph{Ground state
  solutions of {N}ehari-{P}ohozaev type for
  {S}chr\"{o}dinger-{P}oisson-{S}later equation with zero mass and critical
  growth}}.
\newblock J. Geom. Anal., 34(7) Paper No. 221, 19.

\bibitem[{Ianni \& Ruiz(2012)}]{MR2902293}
Ianni I, Ruiz D (2012).
\newblock \href{https://doi.org/10.1142/S0219199712500034}{\emph{Ground and
  bound states for a static {S}chr\"odinger-{P}oisson-{S}later problem}}.
\newblock Commun. Contemp. Math., 14(1) 1250003, 22.

\bibitem[{Kang \emph{et~al.}(2023)Kang, Liu \& Tang}]{Kang2023}
Kang JC, Liu XQ, Tang CL (2023).
\newblock \href{https://doi.org/10.1007/s12220-022-01120-w}{\emph{Ground state
  sign-changing solutions for critical {S}chr\"{o}dinger-{P}oisson system with
  steep potential well}}.
\newblock J. Geom. Anal., 33(2) Paper No. 59, 24.

\bibitem[{Liu \& Perera(2024)}]{MR4785363}
Liu S, Perera K (2024).
\newblock \href{https://doi.org/10.1007/s00526-024-02811-8}{\emph{Multiple
  solutions for $(p,q)$-{L}aplacian equations in $\mathbb{R}^n$ with critical
  or subcritical exponents}}.
\newblock Calc. Var. Partial Differential Equations, 63(8) Paper No. 199, 15.

\bibitem[{Liu \emph{et~al.}(2019)Liu, Zhang \& Huang}]{MR3912770}
Liu Z, Zhang Z, Huang S (2019).
\newblock \href{https://doi.org/10.1016/j.jde.2018.10.048}{\emph{Existence and
  nonexistence of positive solutions for a static
  {S}chr\"odinger-{P}oisson-{S}later equation}}.
\newblock J. Differential Equations, 266(9) 5912--5941.

\bibitem[{Mercuri \emph{et~al.}(2016)Mercuri, Moroz \&
  Van~Schaftingen}]{Mercuri2016}
Mercuri C, Moroz V, Van~Schaftingen J (2016).
\newblock \href{https://doi.org/10.1007/s00526-016-1079-3}{\emph{Groundstates
  and radial solutions to nonlinear {S}chr\"{o}dinger-{P}oisson-{S}later
  equations at the critical frequency}}.
\newblock Calc. Var. Partial Differential Equations, 55(6) Art. 146, 58.

\bibitem[{Mercuri \&
  Perera(2024)}]{mercuri2024variationalmethodsscaledfunctionals}
Mercuri C, Perera K (2024).
\newblock \href{https://arxiv.org/abs/2411.15887}{\emph{Variational methods for
  scaled functionals with applications to the schr\"{o}dinger-poisson-slater
  equation}}.

\bibitem[{Perera(2024)}]{perera2024abstractmultiplicitytheoremsapplications}
Perera K (2024).
\newblock \href{https://arxiv.org/abs/2308.07901}{\emph{Abstract multiplicity
  theorems and applications to critical growth problems}}.

\bibitem[{Ruiz(2010)}]{MR2679375}
Ruiz D (2010).
\newblock \href{http://dx.doi.org/10.1007/s00205-010-0299-5}{\emph{On the
  {S}chr\"odinger-{P}oisson-{S}later system: behavior of minimizers, radial and
  nonradial cases}}.
\newblock Arch. Ration. Mech. Anal., 198(1) 349--368.

\bibitem[{Wang \& Yuan(2023)}]{MR4653880}
Wang Y, Yuan R (2023).
\newblock
  \href{https://doi.org/10.1080/00036811.2022.2130779}{\emph{Nonexistence of
  ground state sign-changing solutions for autonomous {S}chr\"odinger-{P}oisson
  system with critical growth}}.
\newblock Appl. Anal., 102(16) 4652--4658.

\bibitem[{Willem(1996)}]{MR1400007}
Willem M (1996).
\newblock \emph{Minimax theorems}.
\newblock Progress in Nonlinear Differential Equations and their Applications,
  24. Birkh\"auser Boston Inc., Boston, MA.

\bibitem[{Yao \emph{et~al.}(2019)Yao, Li, Zhang \& Mu}]{Yao2019}
Yao X, Li X, Zhang F, Mu C (2019).
\newblock \href{https://doi.org/10.1155/2019/8453176}{\emph{Infinitely many
  solutions of {S}chr\"{o}dinger-{P}oisson equations with critical and
  sublinear terms}}.
\newblock Adv. Math. Phys., Art. ID 8453176, 9.

\bibitem[{Zhang(2025)}]{Zhang2025a}
Zhang J (2025).
\newblock \href{https://doi.org/10.1007/s00526-024-02923-1}{\emph{Multi-bump
  solutions to {S}chr\"{o}dinger-{P}oisson equations with critical growth in
  {${\Bbb {R}}^3$}}}.
\newblock Calc. Var. Partial Differential Equations, 64(2) Paper No. 62, 28.

\bibitem[{Zhao \& Zhao(2009)}]{Zhao2009}
Zhao L, Zhao F (2009).
\newblock \href{https://doi.org/10.1016/j.na.2008.02.116}{\emph{Positive
  solutions for {S}chr\"{o}dinger-{P}oisson equations with a critical
  exponent}}.
\newblock Nonlinear Anal., 70(6) 2150--2164.

\end{thebibliography}
\end{document}